# A TWO-SAMPLE TEST FOR HIGH-DIMENSIONAL DATA WITH APPLICATIONS TO GENE-SET TESTING[1]

By Song Xi Chen and Ying-Li Qin

*Iowa State University and Peking University, and Iowa State University*

We propose a two-sample test for the means of high-dimensional data when the data dimension is much larger than the sample size. Hotelling's classical $T^2$ test does not work for this "large $p$, small $n$" situation. The proposed test does not require explicit conditions in the relationship between the data dimension and sample size. This offers much flexibility in analyzing high-dimensional data. An application of the proposed test is in testing significance for sets of genes which we demonstrate in an empirical study on a leukemia data set.

**1. Introduction.** High-dimensional data are increasingly encountered in many applications of statistics and most prominently in biological and financial studies. A common feature of high-dimensional data is that, while the data dimension is high, the sample size is relatively small. This is the so-called "large $p$, small $n$" phenomenon where $p/n \to \infty$; here $p$ is the data dimension and $n$ is the sample size. The high data dimension ("large $p$") alone has created the need to renovate and rewrite some of the conventional multivariate analysis procedures; these needs only get much greater for "large $p$ small $n$" situations.

A specific "large $p$, small $n$" situation arises when simultaneously testing a large number of hypotheses which is largely motivated by the identification of significant genes in microarray and genetic sequence studies. A natural question is how many hypotheses can be tested simultaneously. This paper tries to answer this question in the context of two-sample simultaneous tests for means. Consider two random samples $X_{i1}, \ldots, X_{in_i} \in R^p$ for $i = 1$ and $2$ which have means $\mu_1 = (\mu_{11}, \ldots, \mu_{1p})^T$ and $\mu_2 = (\mu_{21}, \ldots, \mu_{2p})^T$ and covariance matrices $\Sigma_1$ and $\Sigma_2$, respectively. We consider testing the following

Received October 2008; revised May 2009.
[1]Supported by NSF Grants SES-05-18904, DMS-06-04563 and DMS-07-14978.
*AMS 2000 subject classifications.* Primary 62H15, 60K35; secondary 62G10.
*Key words and phrases.* High dimension, gene-set testing, large $p$ small $n$, martingale central limit theorem, multiple comparison.

high-dimensional hypothesis:

(1.1) $$H_0: \mu_1 = \mu_2 \quad \text{versus} \quad H_1: \mu_1 \neq \mu_2.$$

The hypothesis $H_0$ consists of the $p$ marginal hypotheses $H_{0l}: \mu_{1l} = \mu_{2l}$ for $l = 1, \ldots, p$ regarding the means on each data dimension.

There have been a series of important studies on the high-dimensional problem. Van der Laan and Bryan (2001) show that the sample mean of $p$-dimensional data can consistently estimate the population mean uniformly across $p$ dimensions if $\log(p) = o(n)$ for bounded random variables. In a major generalization, Kosorok and Ma (2007) consider uniform convergence for a range of univariate statistics constructed for each data dimension which includes the marginal empirical distribution, sample mean and sample median. They establish the uniform convergence across $p$ dimensions when $\log(p) = o(n^{1/2})$ or $\log(p) = o(n^{1/3})$, depending on the nature of the marginal statistics. Fan, Hall and Yao (2007) evaluate approximating the overall level of significance for simultaneous testing of means. They demonstrate that the bootstrap can accurately approximate the overall level of significance if $\log(p) = o(n^{1/3})$ when the marginal tests are performed based on the normal or the $t$-distributions. See also Fan, Peng and Huang (2005) and Huang, Wang and Zhang (2005) for high-dimensional estimation and testing in semiparametric regression models.

In an important work, Bai and Saranadasa (1996) propose using $\|\bar{X}_1 - \bar{X}_2\|$ to replace $(\bar{X}_1 - \bar{X}_2)^T S_n^{-1} (\bar{X}_1 - \bar{X}_2)$ in Hotelling's $T^2$-statistic where $\bar{X}_1$ and $\bar{X}_2$ are the two sample means, $S_n$ is the pooled sample covariance by assuming $\Sigma_1 = \Sigma_2 = \Sigma$ and $\|\cdot\|$ denotes the Euclidean norm in $R^p$. They establish the asymptotic normality of the test statistics and show that it has attractive power property when $p/n \to c < \infty$ and under some restriction on the maximum eigenvalue of $\Sigma$. However, the requirement of $p$ and $n$ being of the same order is too restrictive to be used in the "large $p$ small $n$" situation.

To allow simultaneous testing for ultra high-dimensional data, we construct a test which allows $p$ to be arbitrarily large and independent of the sample size as long as, in the case of common covariance $\Sigma$, $\text{tr}(\Sigma^4) = o\{\text{tr}^2(\Sigma^2)\}$ where $\text{tr}(\cdot)$ is the trace operator of a matrix. The above condition on $\Sigma$ is trivially true for any $p$ if either all the eigenvalues of $\Sigma$ are bounded or the largest eigenvalue is of smaller order of $(p-b)^{1/2} b^{-1/4}$ where $b$ is the number of unbounded eigenvalues. We establish the asymptotic normality of a test statistic which leads to a two-sample test for high-dimensional data.

Testing significance for gene-sets rather than a single gene is the latest development in genetic data analysis. A critical need for gene-set testing is to have a multivariate test that is applicable to a wide range of data dimensions (the number of genes in a set). It requires $P$-values for all gene-sets to allow procedures based on either the Bonferroni correction or the false discovery rate [Benjamini and Hochberg (1995)] to take into account the multiplicity in the test. We demonstrate in this paper how to use the proposed test for testing significance for gene-sets. An advantage of the proposed test is in its readily producing $P$-values of significance for each gene-set under study so that the multiplicity of multiple testing can be taken into consideration.

The paper is organized as follows. We outline in Section 2 the framework of the two-sample tests for high-dimensional data and introduce the proposed test statistic. Section 3 provides the theoretical properties of the test. How to apply the proposed test of significance for gene-sets is demonstrated in Section 4 which includes an empirical study on an acute lymphoblastic leukemia data set. Results of simulation studies are reported in Section 5. All the technical details are given in Section 6.

**2. Test statistic.** Suppose we have two independent and identically distributed random samples in $R^p$,

$$\{X_{i1}, X_{i2}, \ldots, X_{in_i}\} \stackrel{\text{i.i.d.}}{\sim} F_i \quad \text{for } i = 1 \text{ and } 2,$$

where $F_i$ is a distribution in $R^p$ with mean $\mu_i$ and covariance $\Sigma_i$. A well-pursued interest in high-dimensional data analysis is to test if the two high-dimensional populations have the same mean or not namely

(2.1) $$H_0 : \mu_1 = \mu_2 \quad \text{vs.} \quad H_1 : \mu_1 \neq \mu_2.$$

The above hypothesis consists of $p$ marginal hypotheses regarding the means of each data dimension. An important question from the point view of multiple testing is how many marginal hypotheses can be tested simultaneously. The works of van der Laan and Bryan (2001), Kosorok and Ma (2007) and Fan, Hall and Yao (2007) are designed to address the question. The existing results show that $p$ can reach the rate of $e^{\alpha n^\beta}$ for some positive constants $\alpha$ and $\beta$. In establishing a rate of the above form, both van der Laan and Bryan (2001) and Kosorok and Ma (2007) assume that the marginal distributions of $F_1$ and $F_2$ are all supported on bounded intervals.

Hotelling's $T^2$ test is the conventional test for the above hypothesis when the dimension $p$ is fixed and is less than $n =: n_1 + n_2 - 2$ and when $\Sigma_1 = \Sigma_2 = \Sigma$, say. Its performance for high-dimensional data is evaluated in Bai and Saranadasa (1996) when $p/n \to c \in [0, 1)$ which reveals a decreasing power as $c$ gets larger. A reason for this negative effect of high-dimension is due to having the inverse of the covariance matrix in the $T^2$ statistic. While standardizing by the covariance brings benefits for data with a fixed dimension, it becomes a liability for high-dimensional data. In particular, the sample covariance matrix $S_n$ may not converge to the population covariance when $p$ and $n$ are of the same order. Indeed, Yin, Bai and Krishnaiah (1988) show that when $p/n \to c$, the smallest and the largest eigenvalues of the sample covariance $S_n$ do not converge to the respective eigenvalues of $\Sigma$. The same phenomenon, but on the weak convergence of the extreme eigenvalues of the sample covariance, is found in Tracy and Widom (1996). When $p > n$, Hotelling's $T^2$ statistic is not defined as $S_n$ may not be invertible.

Our proposed test is motivated by Bai and Saranadasa (1996), who propose testing hypothesis (2.1) under $\Sigma_1 = \Sigma_2 = \Sigma$ based on

$$(2.2) \qquad M_n = (\bar{X}_1 - \bar{X}_2)'(\bar{X}_1 - \bar{X}_2) - \tau \operatorname{tr}(S_n),$$

where $S_n = \frac{1}{n}\sum_{i=1}^{2}\sum_{j=1}^{N_i}(X_{ij} - \bar{X}_i)(X_{ij} - \bar{X}_i)'$ and $\tau = \frac{n_1+n_2}{n_1 n_2}$. The key feature of the Bai and Saranadasa proposal is removing $S_n^{-1}$ in Hotelling's $T^2$ since having $S_n^{-1}$ is no longer beneficial when $p/n \to c > 0$. The subtraction of $\operatorname{tr}(S_n)$ in (2.2) is to make $E(M_n) = \|\mu_1 - \mu_2\|^2$. The asymptotic normality of $M_n$ was established and a test statistic was formulated by standardizing $M_n$ with an estimate of its standard deviation.

The following are the main conditions assumed in Bai–Saranadasa's test:

$$(2.3) \qquad p/n \to c < \infty \quad \text{and} \quad \lambda_p = o(p^{1/2});$$

$$(2.4) \quad n_1/(n_1 + n_2) \to k \in (0,1) \quad \text{and} \quad (\mu_1 - \mu_2)'\Sigma(\mu_1 - \mu_2) = o\{\operatorname{tr}(\Sigma^2)/n\},$$

where $\lambda_p$ denotes the largest eigenvalue of $\Sigma$.

A careful study of the $M_n$ statistic reveals that the restrictions on $p$ and $n$, and on $\lambda_p$ in (2.3) are needed to control terms $\sum_{j=1}^{n_i} X'_{ij} X_{ij}$, $i = 1$ and $2$, in $\|\bar{X}_1 - \bar{X}_2\|^2$. However, these two terms are not useful in the testing. To appreciate this point, let us consider

$$T_n =: \frac{\sum_{i \neq j}^{n_1} X'_{1i} X_{1j}}{n_1(n_1 - 1)} + \frac{\sum_{i \neq j}^{n_2} X'_{2i} X_{2j}}{n_2(n_2 - 1)} - 2\frac{\sum_{i=1}^{n_1} \sum_{j=1}^{n_2} X'_{1i} X_{2j}}{n_1 n_2}$$

after removing $\sum_{j=1}^{n_i} X'_{ij} X_{ij}$ for $i = 1$ and $2$ from $\|\bar{X}_1 - \bar{X}_2\|^2$. Elementary derivations show that

$$E(T_n) = \|\mu_1 - \mu_2\|^2.$$

Hence, $T_n$ is basically all we need for testing. Bai and Saranadasa used $\operatorname{tr}(S_n)$ to offset the two diagonal terms. However, $\operatorname{tr}(S_n)$ itself imposes demands on the dimensionality too.

A derivation in the Appendix shows that under $H_1$ and the condition in (3.4),

$$\operatorname{Var}(T_n) = \left\{ \frac{2}{n_1(n_1-1)} \operatorname{tr}(\Sigma_1^2) + \frac{2}{n_2(n_2-1)} \operatorname{tr}(\Sigma_2^2) + \frac{4}{n_1 n_2} \operatorname{tr}(\Sigma_1 \Sigma_2) \right\}\{1 + o(1)\},$$

where the $o(1)$ term vanishes under $H_0$.

**3. Main results.** We assume, like Bai and Saranadasa (1996), the following general multivariate model:

$$(3.1) \qquad X_{ij} = \Gamma_i Z_{ij} + \mu_i \qquad \text{for } j = 1, \ldots, n_i, \ i = 1 \text{ and } 2,$$

where each $\Gamma_i$ is a $p \times m$ matrix for some $m \geq p$ such that $\Gamma_i \Gamma_i' = \Sigma_i$, and $\{Z_{ij}\}_{j=1}^{n_i}$ are $m$-variate independent and identically distributed (i.i.d.) random vectors satisfying $E(Z_{ij}) = 0$, $\text{Var}(Z_{ij}) = I_m$, the $m \times m$ identity matrix. Furthermore, if we write $Z_{ij} = (z_{ij1}, \ldots, z_{ijm})'$, we assume $E(z_{ijk}^4) = 3 + \Delta < \infty$, and

$$(3.2) \qquad E(z_{ijl_1}^{\alpha_1} z_{ijl_2}^{\alpha_2} \cdots z_{ijl_q}^{\alpha_q}) = E(z_{ijl_1}^{\alpha_1}) E(z_{ijl_2}^{\alpha_2}) \cdots E(z_{ijl_q}^{\alpha_q})$$

for a positive integer $q$ such that $\sum_{l=1}^{q} \alpha_l \leq 8$ and $l_1 \neq l_2 \neq \cdots \neq l_q$. Here $\Delta$ describes the difference between the fourth moments of $z_{ijl}$ and $N(0,1)$. Model (3.1) says that $X_{ij}$ can be expressed as a linear transformation of a $m$-variate $Z_{ij}$ with zero mean and unit variance that satisfies (3.2). Model (3.1) is similar to factor models in multivariate analysis. However, instead of having the number of factors $m < p$ in the conventional multivariate analysis, we require $m \geq p$. This is to allow the basic characteristics of the covariance $\Sigma_i$, for instance its rank and eigenvalues, to not be affected by the transformation. The rank and eigenvalues would be affected if $m < p$. The fact that $m$ is arbitrary offers much flexibility in generating a rich collection of dependence structure. Condition (3.2) means that each $Z_{ij}$ has a kind of pseudo-independence among its components $\{z_{ijl}\}_{l=1}^{m}$. Obviously, if $Z_{ij}$ does have independent components, then (3.2) is trivially true.

We do not assume $\Sigma_1 = \Sigma_2$, as it is a rather strong assumption, and most importantly such an assumption is harder to be verified for high-dimensional data. Testing certain special structures of the covariance matrix when $p$ and $n$ are of the same order have been considered in Ledoit and Wolf (2002) and Schott (2005).

We assume

$$(3.3) \qquad n_1/(n_1 + n_2) \to k \in (0, 1) \qquad \text{as } n \to \infty,$$

$$(3.4) \quad (\mu_1 - \mu_2)' \Sigma_i (\mu_1 - \mu_2) = o[n^{-1} \text{tr}\{(\Sigma_1 + \Sigma_2)^2\}] \qquad \text{for } i = 1 \text{ or } 2,$$

which generalize (2.4) to unequal covariances. Condition (3.4) is obviously satisfied under $H_0$ and implies that the difference between $\mu_1$ and $\mu_2$ is small relative to $n^{-1} \text{tr}\{(\Sigma_1 + \Sigma_2)^2\}$ so that a workable expression for the variance of $T_n$ under $H_0$, and the specified local alternative can be derived. It can be viewed as a high-dimensional version of the local alternative hypotheses. When $p$ is fixed, if we use a standard test for two population means, for instance Hotelling's $T^2$ test, the local alternative hypotheses has the form of $\mu_1 - \mu_2 = \tau n^{-1/2}$ for a nonzero constant vector $\tau \in R^p$. Hotelling's test has nontrivial power under such local alternatives [Anderson (2003)]. If we assume each component of $\mu_1 - \mu_2$ is the same, say $\delta$, then the local alternatives imply $\delta = O(n^{-1/2})$ for a fixed $p$. When the difference is $o(n^{-1/2})$, Hotelling's test has nonpower beyond the level of significance.

To gain insight into (3.4) for high-dimensional situations, let us assume all the eigen-values of $\Sigma_i$ are bounded above from infinity and below away from zero so that $\Sigma_i = I_p$ is a special case of such a regime. Let us also assume, like above, that each component of $\mu_1 - \mu_2$ is the same as a fixed $\delta$, namely $\mu_{1l} - \mu_{2l} = \delta$ for $l = 1, \ldots, p$. Then (3.4) implies $\delta = o(n^{-1/2})$ which is a smaller order than $\delta = O(n^{-1/2})$ for the fixed $p$ case. This can be understood as the high-dimensional data ($p \to \infty$) contain more data information which allows finer resolution in differentiating the two means in each component than that in the fixed $p$ case.

To understand the performance of the test when (3.4) is not valid, we reverse the local alternative condition (3.4) to

$$(3.5) \quad n^{-1}\operatorname{tr}\{(\Sigma_1 + \Sigma_2)^2\} = o\{(\mu_1 - \mu_2)'\Sigma_i(\mu_1 - \mu_2)\} \qquad \text{for } i = 1 \text{ or } 2,$$

implying that the Mahanalobis distance between $\mu_1$ and $\mu_2$ is a larger order than that of $n^{-1}\operatorname{tr}\{(\Sigma_1 + \Sigma_2)^2\}$. This condition can be viewed as a version of fixed alternatives. We will establish asymptotic normally of $T_n$ under either (3.4) or (3.5) in Theorem 1.

The condition we impose on $p$ to replace the first part of (2.3) is

$$(3.6) \quad \operatorname{tr}(\Sigma_i\Sigma_j\Sigma_l\Sigma_h) = o[\operatorname{tr}^2\{(\Sigma_1 + \Sigma_2)^2\}] \qquad \text{for } i, j, l, h = 1 \text{ or } 2,$$

as $p \to \infty$. To appreciate this condition, consider the case of $\Sigma_1 = \Sigma_2 = \Sigma$. Then (3.6) becomes

$$(3.7) \qquad\qquad\qquad \operatorname{tr}(\Sigma^4) = o\{\operatorname{tr}^2(\Sigma^2)\}.$$

Let $\lambda_1 \leq \lambda_2 \leq \cdots \leq \lambda_p$ be the eigenvalues of $\Sigma$. If all eigenvalues are bounded, then (3.7) is trivially true. If, otherwise, there are $b$ unbounded eigenvalues with respect to $p$, and the remaining $p - b$ eigenvalues are bounded above by a finite constant $M$ such that $(p - b) \to \infty$ and $(p - b)\lambda_1^2 \to \infty$, then sufficient conditions for (3.7) are

$$(3.8) \quad \lambda_p = o\{(p-b)^{1/2}\lambda_1 b^{-1/4}\} \quad \text{or} \quad \lambda_p = o\{(p-b)^{1/4}\lambda_1^{1/2}\lambda_{p-b+1}^{1/2}\},$$

where $b$ can be either bounded or diverging to infinity, and the smallest eigen-value $\lambda_1$ can converge to zero. To appreciate these, we note that

$$\frac{\operatorname{tr}(\Sigma^4)}{\operatorname{tr}^2(\Sigma^2)} \leq \frac{(p-b)M^4 + b\lambda_p^4}{(p-b)^2\lambda_1^4 + b^2\lambda_{p-b+1}^4 + 2(p-b)b\lambda_1^2\lambda_{p-b+1}^2}.$$

Hence, the ratio converges to 0 under either condition in (3.8).

The following theorem establishes the asymptotic normality of $T_n$.

THEOREM 1. *Under the assumptions (3.1), (3.2), (3.3), (3.6) and either (3.4) or (3.5),*

$$\frac{T_n - \|\mu_1 - \mu_2\|^2}{\sqrt{\mathrm{Var}(T_n)}} \xrightarrow{d} N(0,1) \qquad \text{as } p \to \infty \text{ and } n \to \infty.$$

The asymptotic normality is attained without imposing any explicit restriction between $p$ and $n$ directly. The only restriction on the dimension is (3.6) or (3.7). As the discussion given just before Theorem 1 suggests, (3.7) is satisfied provided that the number of divergent eigenvalues of $\Sigma$ are not too many, and the divergence is not too fast. The reason for attaining this in the case of high-data-dimension is because the statistic $T_n$ is univariate, despite the fact that the hypothesis $H_0$ is of high dimension. This is different from using a high-dimensional statistic. Indeed, Portnoy (1986) considers the central limit theorem for the $p$-dimensional sample mean $\bar{X}$ and finds that the central limit theorem is not valid if $p$ is not a smaller order of $\sqrt{n}$.

As shown in Section 6.1, $\mathrm{Var}(T_n) = \sigma_n^2 \{1 + o(1)\}$ where, under (3.4),

$$(3.9) \quad \sigma_n^2 =: \sigma_{n1}^2 = \frac{2}{n_1(n_1-1)} \mathrm{tr}(\Sigma_1^2) + \frac{2}{n_2(n_2-1)} \mathrm{tr}(\Sigma_2^2) + \frac{4}{n_1 n_2} \mathrm{tr}(\Sigma_1 \Sigma_2)$$

and under (3.5),

$$(3.10) \quad \sigma_n^2 =: \sigma_{n2}^2 = \frac{4}{n_1}(\mu_1 - \mu_2)'\Sigma_1(\mu_1 - \mu_2) + \frac{4}{n_2}(\mu_1 - \mu_2)'\Sigma_2(\mu_1 - \mu_2).$$

In order to formulate a test procedure based on Theorem 1, $\sigma_{n1}^2$ in (3.9) needs to be estimated. Bai and Saranadasa (1996) used the following estimator for $\mathrm{tr}(\Sigma^2)$ under $\Sigma_1 = \Sigma_2 = \Sigma$:

$$\widehat{\mathrm{tr}(\Sigma^2)} = \frac{n^2}{(n+2)(n-1)} \left\{ \mathrm{tr}\, S_n^2 - \frac{1}{n}(\mathrm{tr}\, S_n)^2 \right\}.$$

Motivated by the benefits of excluding terms like $\sum_{j=1}^{n_i} X'_{ij} X_{ij}$ in the formulation of $T_n$, we propose the following estimator of $\mathrm{tr}(\Sigma_i^2)$ and $\mathrm{tr}(\Sigma_1 \Sigma_2)$:

$$\widehat{\mathrm{tr}(\Sigma_i^2)} = \{n_i(n_i-1)\}^{-1} \mathrm{tr}\left\{ \sum_{j \neq k}^{n_i} (X_{ij} - \bar{X}_{i(j,k)}) X'_{ij} (X_{ik} - \bar{X}_{i(j,k)}) X'_{ik} \right\}$$

and

$$\widehat{\mathrm{tr}(\Sigma_1 \Sigma_2)} = (n_1 n_2)^{-1} \mathrm{tr}\left\{ \sum_{l=1}^{n_1} \sum_{k=1}^{n_2} (X_{1l} - \bar{X}_{1(l)}) X'_{1l} (X_{2k} - \bar{X}_{2(k)}) X'_{2k} \right\},$$

where $\bar{X}_{i(j,k)}$ is the $i$th sample mean after excluding $X_{ij}$ and $X_{ik}$, and $\bar{X}_{i(l)}$ is the $i$th sample mean without $X_{il}$. These are similar to the idea of cross-validation, in that when we construct the deviations of $X_{ij}$ and $X_{ik}$ from the sample mean, both $X_{ij}$ and $X_{ik}$ are excluded from the sample mean calculation. By doing so, the above estimators $\widehat{\mathrm{tr}(\Sigma_i^2)}$ and $\widehat{\mathrm{tr}(\Sigma_1\Sigma_2)}$ can be written as the trace of sums of products of independent matrices. We also note that subtraction of only one sample mean per observation is needed in order to avoid a term like $\|X_{ij}\|^4$ which is harder to control asymptotically without an explicit assumption between $p$ and $n$.

The next theorem shows that the above estimators are ratio-consistent to $\mathrm{tr}(\Sigma_i^2)$ and $\mathrm{tr}(\Sigma_1\Sigma_2)$, respectively.

THEOREM 2. *Under the assumptions (3.1)–(3.4) and (3.6), for $i=1$ or 2,*

$$\frac{\widehat{\mathrm{tr}(\Sigma_i^2)}}{\mathrm{tr}(\Sigma_i^2)} \xrightarrow{p} 1 \quad \text{and} \quad \frac{\widehat{\mathrm{tr}(\Sigma_1\Sigma_2)}}{\mathrm{tr}(\Sigma_1\Sigma_2)} \xrightarrow{p} 1 \quad \text{as } p \text{ and } n \to \infty.$$

A ratio-consistent estimator of $\sigma_{n1}^2$ under $H_0$ is

$$\hat{\sigma}_{n1}^2 = \frac{2}{n_1(n_1-1)}\widehat{\mathrm{tr}(\Sigma_1^2)} + \frac{2}{n_2(n_2-1)}\widehat{\mathrm{tr}(\Sigma_2^2)} + \frac{4}{n_1 n_2}\widehat{\mathrm{tr}(\Sigma_1\Sigma_2)}.$$

This together with Theorem 1 leads to the test statistic,

$$Q_n = T_n/\hat{\sigma}_{n1} \xrightarrow{d} N(0,1) \quad \text{as } p \text{ and } n \to \infty,$$

under $H_0$. The proposed test with an $\alpha$ level of significance rejects $H_0$ if $Q_n > \xi_\alpha$ where $\xi_\alpha$ is the upper $\alpha$ quantile of $N(0,1)$.

Theorems 1 and 2 allow us to discuss the power properties of the proposed test. The discussion is made under (3.4) and (3.5), respectively. The power under the local alternative (3.4) is

$$(3.11) \qquad \beta_{n1}(\|\mu_1-\mu_2\|) = \Phi\left(-\xi_\alpha + \frac{nk(1-k)\|\mu_1-\mu_2\|^2}{\sqrt{2\,\mathrm{tr}\{\tilde{\Sigma}(k)^2\}}}\right),$$

where $\tilde{\Sigma}(k) = (1-k)\Sigma_1 + k\Sigma_2$ and $\Phi$ is the standard normal distribution function. The power of Bai–Saranadasa test has the same form if $\Sigma_1 = \Sigma_2$ and if $p$ and $n$ are of the same order.

The power under (3.5) is

$$\beta_{n2}(\|\mu_1-\mu_2\|) = \Phi\left(-\frac{\sigma_{n1}}{\sigma_{n2}}\xi_\alpha + \frac{\|\mu_1-\mu_2\|^2}{\sigma_{n1}}\right) = \Phi\left(\frac{\|\mu_1-\mu_2\|^2}{\sigma_{n1}}\right)$$

as $\sigma_{n1}/\sigma_{n2} \to 0$. Substitute the expression for $\sigma_{n1}$, and we have

$$(3.12) \qquad \beta_{n2}(\|\mu_1-\mu_2\|) = \Phi\left(\frac{nk(1-k)\|\mu_1-\mu_2\|^2}{\sqrt{2\,\mathrm{tr}\{\tilde{\Sigma}(k)^2\}}}\right).$$

Both (3.11) and (3.12) indicate that the proposed test has nontrivial power under the two cases of the alternative hypothesis as long as

$$n\|\mu_1 - \mu_2\|^2/\sqrt{\mathrm{tr}\{\tilde{\Sigma}(k)^2\}}$$

does not vanish to 0 as $n$ and $p \to \infty$. The flavor of the proposed test is different from tests formulated by combining $p$ marginal tests on $H_{0l}$ [defined after (1.1)] for $l = 1, \ldots, p$. The test statistics of such tests are usually constructed via $\max_{1 \leq l \leq p} T_{nl}$ where $T_{nl}$ is a marginal test statistic for $H_{0l}$. This is the case of Kosorok and Ma (2007) and Fan, Hall and Yao (2007). A condition on $p$ and $n$ is needed to ensure (i) the convergence of $\max_{1 \leq l \leq p} T_{nl}$, and (ii) $p$ can reach an order of $\exp(\alpha n^\beta)$ for positive constants $\alpha$ and $\beta$. Usually some additional assumptions are needed; for instance, Kosorok and Ma (2007) assume each component of the random vector has compact support for testing means.

Naturally, if the number of significant univariate hypotheses ($\mu_{1l} \neq \mu_{2l}$) is a lot less than $p$, which is the so-called sparsity scenario, a simultaneous test like the one we propose may encounter a loss of power. This is actually quantified by the power expression (3.11). Without loss of generality, suppose that each $\mu_i$ can be partitioned as $(\mu_i^{(1)'}, \mu_i^{(2)'})'$ so that under $H_1 : \mu_1^{(1)} = \mu_2^{(1)}$ and $\mu_1^{(2)} \neq \mu_2^{(2)}$ where $\mu_i^{(1)}$ is of $p_1$-dimensional and $\mu_i^{(2)}$ is of $p_2$-dimensional and $p_1 + p_2 = p$. Then $\|\mu_1 - \mu_2\| = p_2 \delta^2$ for some positive constant $\delta^2$. Suppose that $\lambda_{m_0}$ be the smallest nonzero eigenvalue of $\tilde{\Sigma}(k)$. Then under the local alternative (3.4), the asymptotic power is bounded above and below by

$$\Phi\left(-\xi_\alpha + \frac{nk(1-k)p_2\delta^2}{\sqrt{2p}\lambda_p}\right) \leq \beta(\|\mu_1 - \mu_2\|) \leq \Phi\left(-\xi_\alpha + \frac{nk(1-k)p_2\delta^2}{\sqrt{2(p-m_0)}\lambda_{m_0}}\right).$$

If $p$ is very large relative to $n$ and $p_2$ under both high-dimensionality and sparsity, so that $nk(1-k)p_2\eta^2/\sqrt{2(p-m_0)} \to 0$, the test could endure low power. With this in mind, we check on the performance of the test under sparsity in simulation studies in Section 5. The simulations show that the proposed test has a robust power and is in fact more powerful than tests based on multiple comparisons with either the Bonferroni or false discovery rate (FDR) procedures. We note here that, due to the multivariate nature of the test and the hypothesis, the proposed test cannot identify which components are significant after the null multivariate hypothesis is rejected. Additional follow-up procedures have to be employed for that purpose. The proposed test becomes very useful when the purpose is to identify significant groups of components like sets of genes, as illustrated in Section 4. The above discussion can be readily extended to the case of (3.5) due to the similarity in the two power functions.

The proposed two-sample test can be modified for paired observations $\{(Y_{i1}, Y_{i2})\}_{i=1}^n$ where $Y_{i1}$ and $Y_{i2}$ are two measurements of $p$-dimensions on a subject $i$ before and after a treatment. Let $X_i = Y_{i2} - Y_{i1}$, $\mu = E(X_i)$ and $\Sigma = \text{Var}(X_i)$. This is effectively a one-sample problem with high-dimensional data. The hypothesis of interest is

$$H_0 : \mu = 0 \quad \text{vs.} \quad H_1 : \mu \neq 0.$$

We can use $F_n = \sum_{i \neq j}^n X_i' X_j / \{n(n-1)\}$ as the test statistic. It is readily shown that $E(F_n) = \mu'\mu$ and $\text{Var}(F_n) = \frac{2}{n(n-1)} \text{tr}(\Sigma_1^2)\{1 + o(1)\}$ under both $H_0$ and $H_1$ if we assume a condition similar to (3.4) so that $\mu'\Sigma\mu = o\{n^{-1} \text{tr}(\Sigma^2)\}$, and the asymptotic normality of $F_n$ by adding $\text{tr}(\Sigma^4) = o\{\text{tr}^2(\Sigma^2)\}$, a variation of (3.6), can be established by utilizing part of the proof on the asymptotic normality of $T_n$. The $\text{tr}(\Sigma^2)$ can be ratio-consistently estimated with $n_1$ replaced by $n$ in $\widehat{\text{tr}(\Sigma_1^2)}$ which leads to a ratio-consistent variance estimation for $F_n$. Then the test and its power can be written out in similar ways as those for the two-sample test.

When $p = O(1)$, which may be viewed as having finite dimension, the asymptotic normality as conveyed in Theorem 1 may not be valid anymore. It may be shown under conditions (3.1)–(3.4) without (3.6), as condition (3.6) is no longer relevant when $p$ is bounded, that the test statistic $(n_1 + n_2)T_n$ converges to $\sum_{l=1}^{2p} \eta_l \chi_{1,l}^2$ where $\{\chi_{1,l}^2\}_{l=1}^{2p}$ are independent $\chi_1^2$ distributed random variables, and $\{\eta_l\}_{l=1}^{2p}$ is a set of constants. The conclusion of Theorem 2 remains valid when $p$ is bounded. The proposed test can still be used for testing in this situation of bounded dimension with estimated critical values via estimation of $\{\eta_l\}_{l=1}^{2p}$. However, people may like to use a test specially catered for such a case such as, for instance, Hotelling's test.

**4. Gene-set testing.** Identifying sets of genes which are significant with respect to certain treatments is the latest development in genetics research [see Barry, Nobel and Wright (2005), Recknor, Nettleton and Reecy (2008), Efron and Tibshrini (2007) and Newton et al. (2007)]. Biologically speaking, each gene does not function individually in isolation. Rather, one gene tends to work with other genes to achieve certain biological tasks.

Suppose that $S_1, \ldots, S_q$ be $q$ sets of genes, where the gene-set $S_g$ consists of $p_g$ genes. Let $F_{1S_g}$ and $F_{2S_g}$ be the distribution functions corresponding to $S_g$ under the treatment and control, and $\mu_{1S_g}$ and $\mu_{2S_g}$ be their respective means. The hypothesis of interest is

$$H_{0g} : \mu_{1S_g} = \mu_{2S_g} \quad \text{for } g = 1, \ldots, q.$$

The gene sets $\{S_g\}_{g=1}^q$ can overlap as a gene can belong to several functional groups, and $p_g$, the number of genes in a set, can range from a moderate to a very large number. So, there are issues of both multiplicity and high-dimensionality in gene-set testing.

We propose applying the proposed test for the significance of each gene-set $\mathcal{S}_g$ when $p_g$ is large. When $p_g$ is of low-dimension, Hotelling's test may be used. Let $pv_g$, $g = 1,\ldots,q$ be the $P$-values obtained from these tests. To control the overall family-wise error rate, we can employ the Bonferroni procedure; to control FDR, we can use Benjamini and Hochberg's (1995) method or its variations as in Benjamini and Yekutieli (2001) and Storey, Taylor and Siegmund (2004). These lead to control of the family-wise error rate or FDR in the context of gene-sets testing. In contrast, tests based on univariate testing have difficulties in producing $P$-values for gene-sets.

Acute lymphoblastic leukemia (ALL) is a form of leukemia, a cancer of white blood cells. The ALL data [Chiaretti et al. (2004)] contains microarray expressions for 128 patients with either T-cell or B-cell type leukemia. Within the B-cell type leukemia, there are two sub-classes representing two molecular classes: the BCR/ABL class and NEG class. The data set has been analyzed by Dudoit, Keles and van der Laan (2008) using a different technology.

Gene-sets are technically defined in gene ontology (GO) system that provides structured and controlled vocabularies producing names of gene-sets (also called GO terms). There are three groups of gene ontologies of interest: biological processes (BP), cellular components (CC) and molecular functions (MF). We carried out preliminary screening for gene-filtering using the approach in Gentleman et al. (2005), which left 2391 genes for analysis. There are 575 unique GO terms in BP category, 221 in MF and 154 in CC for the ALL data. The largest gene-set contains 2059 genes in BP, 2112 genes in MF and 2078 genes in CC; and the GO terms of the three categories share 1861 common genes. We are interested in detecting differences in the expression levels of gene-sets between the BCR/ABL molecular sub-class ($n_1 = 37$) and the NEG molecular sub-class ($n_2 = 42$) for each of the three categories.

We applied the proposed two-sample test with a 5% significance level to test each of the gene-sets in conjunction with the Bonferroni correction to control the family-wise error rate at 0.05 level. It was found that there were 259 gene-sets declared significant in the BP group, 110 in the MF group and 53 in the CC group. Figure 1 displays the histograms of the $P$-values and the values of test statistic $Q_n$ for the three gene-categories. It shows a strong nonuniform distribution of the $P$-values with a large number of $P$-values clustered near 0. At the same time, the $Q_n$-value plots indicate the average $Q_n$-values are much larger than zero. These explain the large number of significant gene-sets detected by the proposed test.

The number of the differentially expressed gene-sets may seem to be high. This was mainly due to overlapping gene-sets. To appreciate this point, we computed for each (say $i$th) significant gene-set, the number of other significant gene-sets which overlapped with it, say $b_i$; and obtained the average of $\{b_i\}$ and their standard deviation. The average number of overlaps (standard deviation) for BP group was 198.9 (51.3), 55.6 (25.2) for MF and 41.6 (9.5) for CC. These number are indeed very high and reveals the gene-sets and their $P$-values are highly dependent.

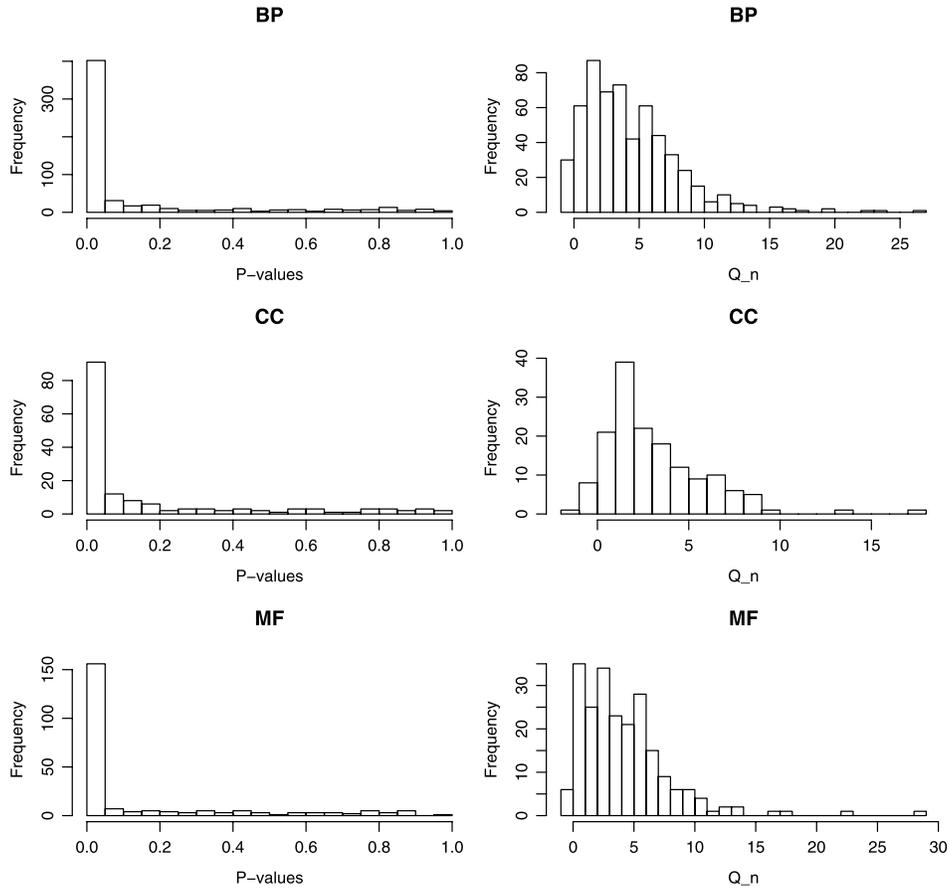

Fig. 1. *Two-sample tests for differentially expressed gene-sets between BCR/ABL and NEG class ALL: histograms of P-values (left panels) and $Q_n$-values (right panels) for BP, CC and MF gene categories.*

Finally, we carried out back-testing for the same hypothesis by randomly splitting the 42 NEG class into two sub-classes of equal sample size and testing for mean differences. This set-up led to the situation of $H_0$. Figure 2 reports the $P$-values and $Q_n$-values for the three gene ontology groups. We note that the distributions of the $P$-values are much closer to the uniform distribution than Figure 1. It is observed that the histograms of $Q_n$-values are centered close to zero and are much closer to the normal distribution than their counterparts in Figure 1 which is reassuring.

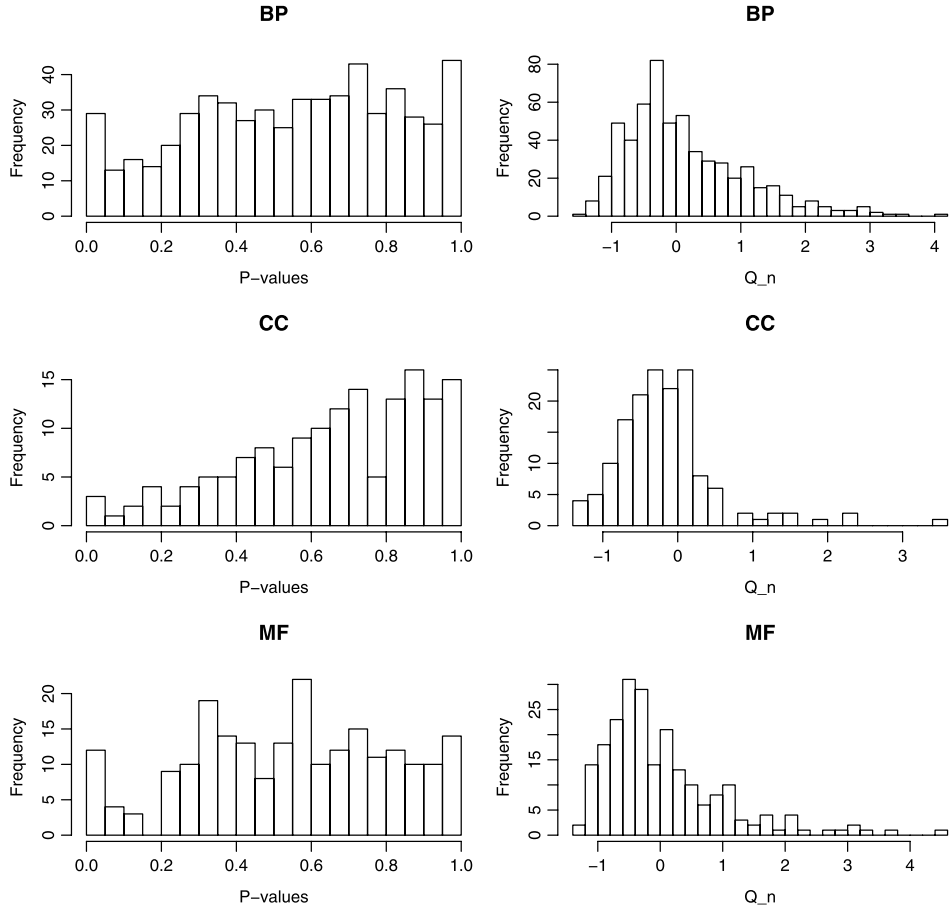

FIG. 2. *Back-testing for differentially expressed gene-sets between two randomly assigned NEG groups: histograms of P-values (left panels) and $Q_n$-values (right panels) for BP, CC and MF gene categories.*

**5. Simulation studies.** In this section, we report results from simulation studies which were designed to evaluate the performance of the proposed two-sample tests for high-dimensional data. For comparison, we also conducted the test proposed by Bai and Saranadasa (1996) (BS test), and two tests based on multiple comparison procedures by employing the Bonferroni and the FDR control [Benjamini and Hochberg (1995)]. The procedure controls the family-wise error rate at a level of significance $\alpha$ which coincides with the significance for the FDR control, the proposed test and the BS test.

In the two multiple comparison procedures, we conducted univariate two-sample $t$-tests for the univariate hypotheses $H_{0l}: \mu_{1l} = \mu_{2l}$ vs $H_{1l}: \mu_{1l} \ne \mu_{2l}$ for $l = 1, 2, \ldots, p$.

Two simulation models for $X_{ij}$ are considered. One has a moving average structure that allows a general dependent structure; the other could allocate the the alternative hypotheses sparsely which enables us to evaluate the performance of the tests under sparsity.

5.1. *Moving average model.* The first simulation model has the following moving average structure:

$$X_{ijk} = \rho_1 Z_{ijk} + \rho_2 Z_{ijk+1} + \cdots + \rho_p Z_{ijk+p-1} + \mu_{ij}$$

for $i = 1$ and $2$, $j = 1, 2, \ldots, n_i$ and $k = 1, 2, \ldots, p$ where $\{Z_{ijk}\}$ are, respectively, i.i.d. random variables. We consider two distributions for the innovations $\{Z_{ijk}\}$. One is a centralized Gamma$(4, 1)$ so that it has zero mean, and the other is $N(0, 1)$.

For each distribution of $\{Z_{ijk}\}$, we consider two configurations of dependence among components of $X_{ij}$. One has weaker dependence with $\rho_l = 0$ for $l > 3$. This prescribes a "two dependence" moving average structure where $X_{ijk_1}$ and $X_{ijk_2}$ are dependent only if $|k_1 - k_2| \le 2$. The $\{\rho_l\}_{l=1}^3$ are generated independently from $U(2, 3)$ which are $\rho_1 = 2.883$, $\rho_2 = 2.794$ and $\rho_3 = 2.849$ and are kept fixed throughout the simulation. The second configuration has all $\rho_l$'s generated from $U(2, 3)$, and again remain fixed throughout the simulation. We call this the "full dependence case." The above dependence structures assigns equal covariance matrices $\Sigma_1 = \Sigma_2 = \Sigma$ and allows a meaningful comparison with the BS test.

Without loss of generality, we fix $\mu_1 = 0$ and choose $\mu_2$ in the same fashion as Benjamini and Hochberg (1995). Specifically, the percentage of true null hypotheses $\mu_{1l} = \mu_{2l}$ for $l = 1, \ldots, p$ were chosen to be 0%, 25%, 50%, 75%, 95% and 99% and 100%, respectively. By experimenting with 95% and 99% we gain information on the performance of the test when $\mu_{1l} \ne \mu_{2l}$ are sparse. It provides empirical checks on the potential concerns of the power of the simultaneous high-dimensional tests as made at the end of Section 3. At each percentage level of true null, three patterns of allocation are considered for the nonzero $\mu_{2l}$ in $\mu_2 = (\mu_{21}, \ldots, \mu_{2p})'$: (i) the equal allocation where all the nonzero $\mu_{2l}$ are equal; (ii) linearly increasing and (iii) linearly decreasing allocations as specified in Benjamini and Hochberg (1995). To make the power comparable among the configurations of $H_1$, we set $\eta =: \|\mu_1 - \mu_2\|^2 / \sqrt{\operatorname{tr}(\Sigma^2)} = 0.1$ throughout the simulation. We chose $p = 500$ and $1000$ and $n = [20 \log(p)] = 124$ and $138$, respectively.

Tables 1 and 2 report the empirical power and size of the four tests with Gamma innovations at a 5% nominal significance level or family-wise error rate or FDR based on 5000 simulations. The results for the normal innovations have a similar pattern, and are not reported here. The simulation results in Tables 1 and 2 can be summarized as follows. The proposed test is much more powerful than the Bai–Saranadasa test for all cases considered in the simulation while maintaining a reasonably-sized approximation to the nominal 5% level. Both the proposed test and the Bai–Saranadasa test are more powerful than the two tests based on the multiple univariate testing using the Bonferroni and FDR procedures. This is expected as both the proposed and Bai–Saranadasa test are designed to test for the entire $p$-dimensional hypotheses while the multiple testing procedures are targeted at the individual univariate hypothesis. What is surprising is that when the percentage of true null is high, at 95% and 99%, the proposed test still is much more powerful than the two multiple testing procedures for all three allocations of the nonzero components in $\mu_2$. It is observed that the sparsity (95% and 99% true null) does reduce the power of the proposed test a little. However, the proposed test still enjoys good power, especially when compared with the other three tests. We also observe that when there is more dependence among multivariate components of the data vectors in the full dependence model, there is a drop in the power for each of the tests. The power of the tests based on the Bonferroni and FDR procedures is alarmingly low and is only slightly larger than the nominal significance level.

Table 1
*Empirical power and size for the 2-dependence model with Gamma innovation*

| Type of allocation | % of true null | $p=500, n=124$ | | | | $p=1000, n=138$ | | | |
|---|---|---|---|---|---|---|---|---|---|
| | | NEW | BS | Bonf | FDR | NEW | BS | Bonf | FDR |
| Equal | 0% | 0.511 | 0.399 | 0.13 | 0.16 | 0.521 | 0.413 | 0.11 | 0.16 |
| | 25% | 0.521 | 0.387 | 0.14 | 0.16 | 0.518 | 0.410 | 0.12 | 0.16 |
| | 50% | 0.513 | 0.401 | 0.13 | 0.17 | 0.531 | 0.422 | 0.12 | 0.17 |
| | 75% | 0.522 | 0.389 | 0.13 | 0.18 | 0.530 | 0.416 | 0.11 | 0.17 |
| | 95% | 0.501 | 0.399 | 0.14 | 0.16 | 0.500 | 0.398 | 0.13 | 0.17 |
| | 99% | 0.499 | 0.388 | 0.13 | 0.15 | 0.507 | 0.408 | 0.15 | 0.18 |
| | 100% (size) | 0.043 | 0.043 | 0.040 | 0.041 | 0.043 | 0.042 | 0.042 | 0.042 |
| Increasing | 0% | 0.520 | 0.425 | 0.11 | 0.13 | 0.522 | 0.409 | 0.12 | 0.15 |
| | 25% | 0.515 | 0.431 | 0.12 | 0.15 | 0.523 | 0.412 | 0.14 | 0.16 |
| | 50% | 0.512 | 0.412 | 0.13 | 0.15 | 0.528 | 0.421 | 0.15 | 0.17 |
| | 75% | 0.522 | 0.409 | 0.15 | 0.17 | 0.531 | 0.431 | 0.16 | 0.19 |
| | 95% | 0.488 | 0.401 | 0.14 | 0.15 | 0.500 | 0.410 | 0.15 | 0.17 |
| | 99% | 0.501 | 0.409 | 0.15 | 0.17 | 0.511 | 0.412 | 0.15 | 0.16 |
| | 100% (size) | 0.042 | 0.041 | 0.040 | 0.041 | 0.042 | 0.040 | 0.039 | 0.041 |
| Decreasing | 0% | 0.522 | 0.395 | 0.11 | 0.15 | 0.533 | 0.406 | 0.09 | 0.15 |
| | 25% | 0.530 | 0.389 | 0.11 | 0.15 | 0.530 | 0.422 | 0.11 | 0.17 |
| | 50% | 0.528 | 0.401 | 0.12 | 0.17 | 0.522 | 0.432 | 0.12 | 0.17 |
| | 75% | 0.533 | 0.399 | 0.13 | 0.18 | 0.519 | 0.421 | 0.12 | 0.17 |
| | 95% | 0.511 | 0.410 | 0.12 | 0.15 | 0.508 | 0.411 | 0.15 | 0.18 |
| | 99% | 0.508 | 0.407 | 0.14 | 0.15 | 0.507 | 0.418 | 0.16 | 0.17 |
| | 100% (size) | 0.041 | 0.042 | 0.041 | 0.042 | 0.042 | 0.040 | 0.040 | 0.042 |

TABLE 2
*Empirical power and size for the full-dependence model with Gamma innovation*

| Type of allocation | % of true null | $p=500, n=124$ | | | | $p=1000, n=138$ | | | |
|---|---|---|---|---|---|---|---|---|---|
| | | NEW | BS | Bonf | FDR | NEW | BS | Bonf | FDR |
| Equal | 0% | 0.322 | 0.120 | 0.08 | 0.10 | 0.402 | 0.216 | 0.09 | 0.11 |
| | 25% | 0.318 | 0.117 | 0.08 | 0.10 | 0.400 | 0.218 | 0.08 | 0.11 |
| | 50% | 0.316 | 0.115 | 0.09 | 0.11 | 0.409 | 0.221 | 0.09 | 0.10 |
| | 75% | 0.307 | 0.113 | 0.10 | 0.12 | 0.410 | 0.213 | 0.09 | 0.13 |
| | 95% | 0.233 | 0.128 | 0.11 | 0.14 | 0.308 | 0.215 | 0.10 | 0.13 |
| | 99% | 0.225 | 0.138 | 0.12 | 0.15 | 0.316 | 0.207 | 0.11 | 0.12 |
| | 100% (size) | 0.041 | 0.041 | 0.043 | 0.043 | 0.042 | 0.042 | 0.040 | 0.041 |
| Increasing | 0% | 0.331 | 0.121 | 0.09 | 0.12 | 0.430 | 0.225 | 0.10 | 0.11 |
| | 25% | 0.336 | 0.119 | 0.10 | 0.12 | 0.423 | 0.231 | 0.12 | 0.12 |
| | 50% | 0.329 | 0.123 | 0.12 | 0.14 | 0.422 | 0.226 | 0.13 | 0.14 |
| | 75% | 0.330 | 0.115 | 0.12 | 0.15 | 0.431 | 0.222 | 0.14 | 0.15 |
| | 95% | 0.219 | 0.120 | 0.12 | 0.13 | 0.311 | 0.218 | 0.14 | 0.15 |
| | 99% | 0.228 | 0.117 | 0.13 | 0.15 | 0.315 | 0.217 | 0.15 | 0.17 |
| | 100% (size) | 0.041 | 0.040 | 0.042 | 0.043 | 0.042 | 0.042 | 0.040 | 0.042 |
| Decreasing | 0% | 0.320 | 0.117 | 0.08 | 0.11 | 0.411 | 0.213 | 0.08 | 0.10 |
| | 25% | 0.323 | 0.119 | 0.09 | 0.11 | 0.408 | 0.210 | 0.08 | 0.11 |
| | 50% | 0.327 | 0.120 | 0.11 | 0.12 | 0.403 | 0.208 | 0.09 | 0.10 |
| | 75% | 0.322 | 0.122 | 0.12 | 0.12 | 0.400 | 0.211 | 0.12 | 0.13 |
| | 95% | 0.217 | 0.109 | 0.12 | 0.15 | 0.319 | 0.207 | 0.12 | 0.15 |
| | 99% | 0.224 | 0.111 | 0.13 | 0.16 | 0.327 | 0.205 | 0.11 | 0.13 |
| | 100% (size) | 0.042 | 0.043 | 0.039 | 0.041 | 0.042 | 0.211 | 0.040 | 0.041 |

We also collected information on the quality of $\text{tr}(\Sigma^2)$ estimation. Table 3 reports empirical averages and standard deviation of $\widehat{\text{tr}(\Sigma^2)}/\text{tr}(\Sigma^2)$. It shows that the proposed estimator for $\text{tr}(\Sigma^2)$ has a much smaller bias and standard deviation than those proposed in Bai and Saranadasa (1996) in all cases, and provides an empirical verification for Theorem 2.

5.2. *Sparse model.* An examination of the previous simulation setting reveals that the strength of the "signals" $\mu_{2l} - \mu_{1l}$ corresponding to the alternative hypotheses are low relative to the level of noise (variance) which may not be a favorable situation for the two tests based on multiple univariate testing. To gain more information on the performance of the tests under sparsity, we consider the following simulation model such that

$$X_{1il} = Z_{1il} \quad \text{and} \quad X_{2il} = \mu_l + Z_{2il} \quad \text{for } l = 1, \ldots, p,$$

TABLE 3
*Empirical averages of $\widehat{\mathrm{tr}(\Sigma^2)}/\mathrm{tr}(\Sigma^2)$ with standard deviations in the parentheses*

| Type of innovation | Type of dependence | $p=500, n=124$ | | |
|---|---|---|---|---|
| | | **NEW** | **BS** | $\mathrm{tr}(\Sigma^2)$ |
| Normal | 2-dependence | 1.03 (0.015) | 1.39 (0.016) | 3102 |
| | Full-dependence | 1.008 (0.00279) | 1.17 (0.0032) | 35,911 |
| Gamma | 2-dependence | 1.03 (0.006) | 1.10 (0.007) | 14,227 |
| | Full-dependence | 1.108 (0.0019) | 1.248 (0.0017) | 152,248 |
| | | $p=1000, n=138$ | | |
| Normal | 2-dependence | 0.986 (0.0138) | 1.253 (0.0136) | 6563 |
| | Full-dependence | 0.995 (0.0026) | 1.072 (0.0033) | 76,563 |
| Gamma | 2-dependence | 1.048 (0.005) | 1.138 (0.006) | 32,104 |
| | Full-dependence | 1.088 (0.00097) | 1.231 (0.0013) | 325,879 |

where $\{Z_{1il}, Z_{2il}\}_{l=1}^{p}$ are mutually independent $N(0,1)$ random variables, and the "signals,"

$$\mu_l = \varepsilon\sqrt{2\log(p)} \quad \text{for } l=1,\ldots,q=[p^c] \quad \text{and} \quad \mu_l = 0 \quad \text{for } l > q,$$

for some $c \in (0,1)$. Here $q$ is the number of significant alternative hypotheses. The sparsity of the hypotheses is determined by $c$: the smaller the $c$ is, the more sparse the alternative hypotheses with $\mu_l \neq 0$. This simulation model is similar to the one used in Abramovich et al. (2006).

According to (3.11), the power of the proposed test has the asymptotic power

$$\beta(\|\mu\|) = \Phi\left(-\xi_\alpha + \frac{np^{(c-1/2)}\varepsilon^2 \log(p)}{2\sqrt{2}}\right),$$

which indicates that the test has a much reduced power if $c < 1/2$ with respect to $p$. We, therefore, chose $p = 1000$ and $c = 0.25, 0.35, 0.45$ and $0.55$, respectively, which leads to $q = 6, 11, 22$, and $44$, respectively. We call $c = 0.25, 0.35$ and $0.45$ the sparse cases.

In order to prevent trivial powers of $\alpha$ or 1 in the simulation, we set $\varepsilon = 0.25$ for $c = 0.25$ and $0.45$; and $\varepsilon = 0.15$ for $c = 0.35$ and $0.55$. Table 4 summarizes the simulations results based on 500 simulations. It shows that in the extreme sparse cases of $c = 0.25$, the FDR and Bonferroni tests have lower power than the proposed test. The power is largely similar among the three tests for $c = 0.35$. However, when the sparsity is moderated to $c = 0.45$, the proposed test starts to surpass the FDR and Bonferroni procedures. The gap in power performance is further increased when $c = 0.55$. Table 5 reports the quality of the variance estimation in Table 5 which shows the proposed variance estimators incur very little bias and variance for even very small sample sizes of $n_1 = n_2 = 10$.

TABLE 4
*Empirical power and size for the sparse model*

| Sample size $(n_1 = n_2)$ | Methods | $\varepsilon = 0.25$ | | | | $\varepsilon = 0.15$ | | | |
|---|---|---|---|---|---|---|---|---|---|
| | | $c = 0.25$ | | $c = 0.45$ | | $c = 0.35$ | | $c = 0.55$ | |
| | | Power | Size | Power | Size | Power | Size | Power | Size |
| 10 | FDR  | 0.084 | 0.056 | 0.180 | 0.040 | 0.044 | 0.034 | 0.066 | 0.034 |
|    | Bonf | 0.084 | 0.056 | 0.170 | 0.040 | 0.044 | 0.034 | 0.062 | 0.032 |
|    | New  | 0.100 | 0.046 | 0.546 | 0.056 | 0.072 | 0.064 | 0.344 | 0.064 |
| 20 | FDR  | 0.380 | 0.042 | 0.855 | 0.044 | 0.096 | 0.036 | 0.326 | 0.058 |
|    | Bonf | 0.368 | 0.038 | 0.806 | 0.044 | 0.092 | 0.034 | 0.308 | 0.056 |
|    | New  | 0.238 | 0.052 | 0.976 | 0.042 | 0.106 | 0.052 | 0.852 | 0.046 |
| 30 | FDR  | 0.864 | 0.042 | 1     | 0.060 | 0.236 | 0.048 | 0.710 | 0.038 |
|    | Bonfe| 0.842 | 0.038 | 0.996 | 0.060 | 0.232 | 0.048 | 0.660 | 0.038 |
|    | New  | 0.408 | 0.050 | 0.998 | 0.058 | 0.220 | 0.054 | 0.988 | 0.042 |

## 6. Technical details.

6.1. *Derivations for $E(T_n)$ and $\mathrm{Var}(T_n)$.* As

$$T_n = \frac{\sum_{i \neq j}^{n_1} X'_{1i} X_{1j}}{n_1(n_1 - 1)} + \frac{\sum_{i \neq j}^{n_2} X'_{2i} X_{2j}}{n_2(n_2 - 1)} - 2 \frac{\sum_{i=1}^{n_1} \sum_{j=1}^{n_2} X'_{1i} X_{2j}}{n_1 n_2},$$

it is straightforward to show that $E(T_n) = \mu'_1 \mu_1 + \mu'_2 \mu_2 - 2\mu'_1 \mu_2 = \|\mu_1 - \mu_2\|^2$.

Let $P_1 = \frac{\sum_{i \neq j}^{n_1} X'_{1i} X_{1j}}{n_1(n_1-1)}$, $P_2 = \frac{\sum_{i \neq j}^{n_2} X'_{2i} X_{2j}}{n_2(n_2-1)}$ and $P_3 = -2 \frac{\sum_{i=1}^{n_1} \sum_{j=1}^{n_2} X'_{1i} X_{2j}}{n_1 n_2}$. It can be shown that

$$\mathrm{Var}(P_1) = \frac{2}{n_1(n_1 - 1)} \mathrm{tr}(\Sigma_1^2) + \frac{4\mu'_1 \Sigma_1 \mu_1}{n_1},$$

TABLE 5
*Average ratios of $\widehat{\sigma_M^2}/\sigma_M^2$ and their standard deviation (in parenthesis) for the sparse model*

| Sample size | True $\sigma_M^2$ | $\varepsilon = 0.25$ | | $\varepsilon = 0.15$ | |
|---|---|---|---|---|---|
| | | $c = 0.25$ | $c = 0.45$ | $c = 0.35$ | $c = 0.55$ |
| $n_1 = n_2 = 10$ | 84.4 | 1.003 (0.0123) | 1.005 (0.0116) | 0.998 (0.0120) | 0.999 (0.0110) |
| $n_1 = n_2 = 20$ | 20.5 | 1.003 (0.0033) | 1.000 (0.0028) | 1.003 (0.0028) | 1.002 (0.0029) |
| $n_1 = n_2 = 30$ | 9.0  | 0.996 (0.0013) | 0.998 (0.0013) | 1.004 (0.0014) | 0.999 (0.0013) |

A TWO-SAMPLE TEST FOR HIGH-DIMENSIONAL DATA 19$$\mathrm{Var}(P_2) = \frac{2}{n_2(n_2-1)} \mathrm{tr}(\Sigma_2^2) + \frac{4\mu_2'\Sigma_2\mu_2}{n_2}$$

and

$$\mathrm{Var}(P_3) = \frac{4}{n_1 n_2} \mathrm{tr}(\Sigma_1 \Sigma_2) + \frac{4\mu_2'\Sigma_1\mu_2}{n_1} + \frac{4\mu_1'\Sigma_2\mu_1}{n_2}.$$

Because the two samples are independent, $\mathrm{Cov}(P_1, P_2) = 0$. Also,

$$\mathrm{Cov}(P_1, P_3) = -\frac{4\mu_1'\Sigma_1\mu_2}{n_1} \quad \text{and} \quad \mathrm{Cov}(P_2, P_3) = -\frac{4\mu_1'\Sigma_2\mu_2}{n_2}.$$

In summary,

$$\mathrm{Var}(T_n) = \frac{2}{n_1(n_1-1)} \mathrm{tr}(\Sigma_1^2) + \frac{2}{n_2(n_2-1)} \mathrm{tr}(\Sigma_2^2) + \frac{4}{n_1 n_2} \mathrm{tr}(\Sigma_1 \Sigma_2)$$
$$+ \frac{4}{n_1}(\mu_1 - \mu_2)'\Sigma_1(\mu_1 - \mu_2) + \frac{4}{n_2}(\mu_1 - \mu_2)'\Sigma_2(\mu_1 - \mu_2).$$

Thus, under $H_0$,

$$\mathrm{Var}(T_n) = \sigma_{n1}^2 =: \frac{2}{n_1(n_1-1)} \mathrm{tr}(\Sigma_1^2) + \frac{2}{n_2(n_2-1)} \mathrm{tr}(\Sigma_2^2) + \frac{4}{n_1 n_2} \mathrm{tr}(\Sigma_1 \Sigma_2).$$

Under $H_1 : \mu_1 \neq \mu_2$, with (3.4),

$$\mathrm{Var}(T_n) = \sigma_{n1}^2 \{1 + o(1)\};$$

and with (3.5),

$$\mathrm{Var}(T_n) = \sigma_{n2}^2 \{1 + o(1)\},$$

where $\sigma_{n2} = \frac{4}{n_1}(\mu_1 - \mu_2)'\Sigma_1(\mu_1 - \mu_2) + \frac{4}{n_2}(\mu_1 - \mu_2)'\Sigma_2(\mu_1 - \mu_2)$.

6.2. *Asymptotic normality of* $T_n$. We note that $T_n = T_{n1} + T_{n2}$ where

(6.1)
$$T_{n1} = \frac{\sum_{i \neq j}^{n_1}(X_{1i} - \mu_1)'(X_{1j} - \mu_1)}{n_1(n_1 - 1)} + \frac{\sum_{i \neq j}^{n_2}(X_{2i} - \mu_2)'(X_{2j} - \mu_2)}{n_2(n_2 - 1)}$$
$$- 2\frac{\sum_{i=1}^{n_1}\sum_{j=1}^{n_2}(X_{1i} - \mu_1)'(X_{2j} - \mu_2)}{n_1 n_2}$$

and

$$T_{n2} = \frac{2\sum_{i=1}^{n_1}(X_{1i} - \mu_1)'(\mu_1 - \mu_2)}{n_1} + \frac{2\sum_{i=1}^{n_2}(X_{2i} - \mu_2)'(\mu_2 - \mu_1)}{n_2}$$
$$+ \mu_1'\mu_1 + \mu_2'\mu_2 - 2\mu_1'\mu_2.$$

It is easy to show that $E(T_{n1}) = 0$ and $E(T_{n2}) = \|\mu_1 - \mu_2\|^2$, and

$$\mathrm{Var}(T_{n2}) = 4n_1^{-1}(\mu_1 - \mu_2)'\Sigma_1(\mu_1 - \mu_2) + 4n_2^{-1}(\mu_2 - \mu_1)'\Sigma_2(\mu_2 - \mu_1).$$

Under (3.4), as

(6.2)
$$\text{Var}\left(\frac{T_{n2} - \|\mu_1 - \mu_2\|^2}{\sigma_{n1}}\right) = o(1),$$

$$\frac{T_n - \|\mu_1 - \mu_2\|^2}{\sqrt{\text{Var}(T_n)}} = \frac{T_{n1}}{\sigma_{n1}} + o_p(1).$$

Under (3.5),

(6.3)
$$\frac{T_n - \|\mu_1 - \mu_2\|^2}{\sqrt{\text{Var}(T_n)}} = \frac{T_{n2} - \|\mu_1 - \mu_2\|^2}{\sigma_{n2}} + o_p(1).$$

As $T_{n2}$ are independent sample averages, its asymptotic normality is readily attainable as shown later. The main task of the proof is for the case under (3.4) when $T_{n1}$ is the contributor of the asymptotic distribution. From (6.1), in the derivation for the asymptotic normality of $T_{n1}$, we can assume without loss of generality that $\mu_1 = \mu_2 = 0$.

Let $Y_i = X_{1i}$ for $i = 1, \ldots, n_1$ and $Y_{j+n_1} = X_{2j}$ for $j = 1, \ldots, n_2$, and for $i \neq j$

$$\phi_{ij} = \begin{cases} n_1^{-1}(n_1 - 1)^{-1} Y_i' Y_j, & \text{if } i, j \in \{1, 2, \ldots, n_1\}, \\ -n_1^{-1} n_2^{-1} Y_i' Y_j, & \text{if } i \in \{1, 2, \ldots, n_1\} \\ & \text{and } j \in \{n_1 + 1, \ldots, n_1 + n_2\}, \\ n_2^{-1}(n_2 - 1)^{-1} Y_i' Y_j, & \text{if } i, j \in \{n_1 + 1, \ldots, n_1 + n_2\}. \end{cases}$$

Define $V_{nj} = \sum_{i=1}^{j-1} \phi_{ij}$ for $j = 2, 3, \ldots, n_1 + n_2$, $S_{nm} = \sum_{j=2}^{m} V_{nj}$ and $\mathcal{F}_{nm} = \sigma\{Y_1, Y_2, \ldots, Y_m\}$ which is the $\sigma$ algebra generated by $\{Y_1, Y_2, \ldots, Y_m\}$. Now

$$T_n = 2 \sum_{j=2}^{n_1+n_2} V_{nj}.$$

LEMMA 1. *For each $n$, $\{S_{nm}, \mathcal{F}_{nm}\}_{m=1}^{n}$ is the sequence of zero mean and a square integrable martingale.*

PROOF. It's obvious that $\mathcal{F}_{nj-1} \subseteq \mathcal{F}_{nj}$, for any $1 \leq j \leq n$ and $S_{nm}$ is of zero mean and square integrable. We only need to show $E(S_{nq}|\mathcal{F}_{nm}) = S_{nm}$ for any $q \geq m$. We note that if $j \leq m \leq n$, then $E(V_{nj}|\mathcal{F}_{nm}) = \sum_{i=1}^{j-1} E(\phi_{ij}|\mathcal{F}_{nm}) = \sum_{i=1}^{j-1} \phi_{ij} = V_{nj}$. If $j > m$, then $E(\phi_{ij}|\mathcal{F}_{nm}) = E(Y_i' Y_j|\mathcal{F}_{nm})$.

If $i > m$, as $Y_i$ and $Y_j$ are both independent of $\mathcal{F}_{nm}$,

$$E(\phi_{ij}|\mathcal{F}_{nm}) = E(\phi_{ij}) = 0.$$

If $i \leq m$, $E(\phi_{ij}|\mathcal{F}_{n,m}) = E(Y_i' Y_j|\mathcal{F}_{n,m}) = Y_i' E(Y_j) = 0$. Hence,

$$E(V_{nj}|\mathcal{F}_{n,m}) = 0.$$

In summary, for $q > m$, $E(S_{nq}|\mathcal{F}_{nm}) = \sum_{j=1}^{q} E(V_{nj}|\mathcal{F}_{nm}) = \sum_{j=1}^{m} V_{nj} = S_{nm}$. This completes the proof of the lemma. $\square$

LEMMA 2. *Under condition (3.4),*
$$\frac{\sum_{j=2}^{n_1+n_2} E[V_{nj}^2|\mathcal{F}_{n,j-1}]}{\sigma_{n_1}^2} \xrightarrow{P} \frac{1}{4}.$$

PROOF. Note that
$$E(V_{nj}^2|\mathcal{F}_{nj-1}) = E\left\{\left(\sum_{i=1}^{j-1} Y_i' Y_j\right)^2 \bigg| \mathcal{F}_{nj-1}\right\} = E\left(\sum_{i_1,i_2=1}^{j-1} Y_{i_1}' Y_j Y_j' Y_{i_2} \bigg| \mathcal{F}_{nj-1}\right)$$
$$= \sum_{i_1,i_2=1}^{j-1} Y_{i_1}' E(Y_j Y_j'|\mathcal{F}_{nj-1}) Y_{i_2} = \sum_{i_1,i_2=1}^{j-1} Y_{i_1}' E(Y_j Y_j') Y_{i_2}$$
$$= \sum_{i_1,i_2=1}^{j-1} Y_{i_1}' \frac{\tilde{\Sigma}_j}{\tilde{n}_j(\tilde{n}_j - 1)} Y_{i_2},$$

where $\tilde{\Sigma}_j = \Sigma_1, \tilde{n}_j = n_1$, for $j \in [1, n_1]$ and $\tilde{\Sigma}_j = \Sigma_2, \tilde{n}_j = n_2$, if $j \in [n_1 + 1, n_1 + n_2]$.

Define
$$\eta_n = \sum_{j=2}^{n_1+n_2} E(V_{nj}^2|\mathcal{F}_{nj-1}).$$

Then
$$E(\eta_n) = \frac{\text{tr}(\Sigma_1^2)}{2n_1(n_1 - 1)} + \frac{\text{tr}(\Sigma_2^2)}{2n_2(n_2 - 1)} + \frac{\text{tr}(\Sigma_1 \Sigma_2)}{(n_1 - 1)(n_2 - 1)}$$
(6.4)
$$= \frac{1}{4}\sigma_{n_1}^2\{1 + o(1)\}.$$

Now consider
$$E(\eta_n^2) = E\left\{\sum_{j=2}^{n_1+n_2} \sum_{i_1,i_2=1}^{j-1} Y_{i_1}' \frac{\tilde{\Sigma}_j}{\tilde{n}_j(\tilde{n}_j - 1)} Y_{i_2}\right\}^2$$
$$= E\left\{2 \sum_{2 \le j_1 < j_2}^{n_1+n_2} \sum_{i_1,i_2=1}^{j_1-1} \sum_{i_3,i_4=1}^{j_2-1} Y_{i_1}' \frac{\tilde{\Sigma}_{j_1}}{\tilde{n}_{j_1}(\tilde{n}_{j_1} - 1)} Y_{i_2} Y_{i_3}' \frac{\tilde{\Sigma}_{j_2}}{\tilde{n}_{j_2}(\tilde{n}_{j_2} - 1)} Y_{i_4}\right.$$
(6.5)
$$\left. + \sum_{j=2}^{n_1+n_2} \sum_{i_1,i_2=1}^{j-1} \sum_{i_3,i_4=1}^{j-1} Y_{i_1}' \frac{\tilde{\Sigma}_j}{\tilde{n}_j(\tilde{n}_j - 1)} Y_{i_2} Y_{i_3}' \frac{\tilde{\Sigma}_j}{\tilde{n}_j(\tilde{n}_j - 1)} Y_{i_4}\right\}$$
$$= 2E(A) + E(B), \quad \text{say,}$$

where

$$A = \sum_{2 \le j_1 < j_2}^{n_1+n_2} \sum_{i_1,i_2=1}^{j_1-1} \sum_{i_3,i_4=1}^{j_2-1} Y'_{i_1} \frac{\tilde{\Sigma}_{j_1}}{\tilde{n}_{j_1}(\tilde{n}_{j_1}-1)} Y_{i_2} Y'_{i_3} \frac{\tilde{\Sigma}_{j_2}}{\tilde{n}_{j_2}(\tilde{n}_{j_2}-1)} Y_{i_4},$$

(6.6)

$$B = \sum_{j=2}^{n_1+n_2} \sum_{i_1,i_2=1}^{j-1} \sum_{i_3,i_4=1}^{j-1} Y'_{i_1} \frac{\tilde{\Sigma}_j}{\tilde{n}_j(\tilde{n}_j-1)} Y_{i_2} Y'_{i_3} \frac{\tilde{\Sigma}_j}{\tilde{n}_j(\tilde{n}_j-1)} Y_{i_4}.$$

Derivations given in Chen and Qin (2008) show

$$2E(A) = \left\{ \frac{\mathrm{tr}^2(\Sigma_1^2)}{4n_1^2(n_1-1)^2} + \frac{\mathrm{tr}^2(\Sigma_2^2)}{4n_2^2(n_2-1)^2} + \frac{\mathrm{tr}(\Sigma_1^2)\mathrm{tr}(\Sigma_1\Sigma_2)}{n_1^2(n_1-1)(n_2-1)} \right.$$
$$+ \frac{\mathrm{tr}(\Sigma_2^2)\mathrm{tr}(\Sigma_1\Sigma_2)}{(n_1-1)n_2(n_2-1)} + \frac{\mathrm{tr}^2(\Sigma_2\Sigma_1)}{n_1 n_2(n_1-1)(n_2-1)}$$
$$\left. + \frac{\mathrm{tr}(\Sigma_1^2)\mathrm{tr}(\Sigma_2^2)}{2n_1(n_1-1)n_2(n_2-1)} \right\} \{1+o(1)\},$$

and $E(B) = o(\sigma_{n_1}^2)$. Hence, from (6.5) and (6.6),

(6.7)
$$E(\eta_n^2) = \left\{ \frac{\mathrm{tr}^2(\Sigma_1^2)}{4n_1^2(n_1-1)^2} + \frac{\mathrm{tr}^2(\Sigma_2^2)}{4n_2^2(n_2-1)^2} + \frac{\mathrm{tr}(\Sigma_1^2)\mathrm{tr}(\Sigma_1\Sigma_2)}{n_1^2(n_1-1)(n_2-1)} \right.$$
$$+ \frac{\mathrm{tr}(\Sigma_2^2)\mathrm{tr}(\Sigma_1\Sigma_2)}{(n_1-1)n_2(n_2-1)} + \frac{\mathrm{tr}^2(\Sigma_2\Sigma_1)}{n_1 n_2(n_1-1)(n_2-1)}$$
$$\left. + \frac{\mathrm{tr}(\Sigma_1^2)\mathrm{tr}(\Sigma_2^2)}{2n_1(n_1-1)n_2(n_2-1)} \right\} + o(\sigma_{n_1}^4).$$

Based on (6.4) and (6.7),

(6.8) $$\mathrm{Var}(\eta_n) = E(\eta_n^2) - E^2(\eta_n) = o(\sigma_{n_1}^4).$$

Combine (6.4) and (6.8), and we have

$$\sigma_{n_1}^{-2} E\left\{ \sum_{j=1}^{n_1+n_2} E(V_{nj}^2 | \mathcal{F}_{n,j-1}) \right\} = \sigma_{n_1}^{-2} E(\eta_n) = \frac{1}{4}$$

and

$$\sigma_{n_1}^{-4} \mathrm{Var}\left\{ \sum_{j=1}^{n_1+n_2} E(V_{nj}^2 | \mathcal{F}_{n,j-1}) \right\} = \sigma_{n_1}^{-4} \mathrm{Var}(\eta_n) = o(1).$$

This completes the proof of Lemma 2.  □

LEMMA 3. *Under condition (3.4),*

$$\sum_{j=2}^{n_1+n_2} \sigma_{n_1}^{-2} E\{V_{nj}^2 I(|V_{nj}| > \epsilon\sigma_{n_1})|F_{nj-1}\} \xrightarrow{p} 0.$$

PROOF. We note that

$$\sum_{j=2}^{n_1+n_2} \sigma_{n_1}^{-2} E\{V_{nj}^2 I(|V_{nj}| > \epsilon\sigma_{n_1})|F_{nj-1}\} \leq \sigma_{n_1}^{-q} \epsilon^{2-q} \sum_{j=1}^{n_1+n_2} E(V_{nj}^q|F_{nj-1}),$$

for some $q > 2$. By choosing $q = 4$, the conclusion of the lemma is true if we can show

(6.9) $$E\left\{\sum_{j=2}^{n_1+n_2} E(V_{nj}^4|F_{nj-1})\right\} = o(\sigma_{n_1}^4).$$

We notice that

$$E\left\{\sum_{j=2}^{n_1+n_2} E(V_{nj}^4|F_{nj-1})\right\} = \sum_{j=2}^{n_1+n_2} E(V_{nj}^4) = O(n^{-8}) \sum_{j=2}^{n_1+n_2} E\left(\sum_{i=1}^{j-1} \phi_{ij}\right)^4$$

The last term can be decomposed as $3Q + P$ where

$$Q = O(n^{-8}) \sum_{j=2}^{n_1+n_2} \sum_{s \neq t}^{j-1} E(Y_j' Y_s Y_s' Y_j Y_j' Y_t Y_t' Y_j)$$

and $P = O(n^{-8}) \sum_{j=2}^{n_1+n_2} \sum_{s=1}^{j-1} E(Y_s' Y_j)^4$. Now (6.9) is true if $3Q + P = o(\sigma_{n_1}^4)$.

Note that

$$Q = O(n^{-8}) \sum_{j=2}^{n_1+n_2} \sum_{s \neq t}^{j-1} E\{\text{tr}(Y_j Y_j' Y_t Y_t' Y_j Y_j' Y_s Y_s')\}$$

$$= O(n^{-4})\left\{\sum_{j=2}^{n_1} \sum_{s \neq t}^{j-1} E(Y_j' \Sigma_1 Y_j Y_j' \Sigma_1 Y_j) + \sum_{j=n_1+1}^{n_1+n_2} \sum_{s \neq t}^{j-1} E(Y_j' \Sigma_t Y_j Y_j' \Sigma_s Y_j)\right\}$$

$$= o(\sigma_{n_1}^4).$$

The last equation follows the similar procedure in Lemma 2 under (3.4).

It remains to show that $P = O(n^{-8}) \sum_{j=2}^{n_1+n_2} \sum_{s=1}^{j-1} E(Y_s' Y_j)^4 = o(\sigma_{n_1}^4)$. Note that

$$P = O(n^{-8}) \sum_{j=2}^{n_1+n_2} \sum_{s=1}^{j-1} E(Y_s' Y_j)^4$$

$$= O(n^{-8}) \sum_{j=2}^{n_1} \sum_{s=1}^{j-1} E(Y_s'Y_j)^4 + O(n^{-8}) \sum_{j=n_1+1}^{n_1+n_2} \sum_{s=1}^{j-1} E(Y_s'Y_j)^4$$

$$= O(n^{-8}) \left\{ \sum_{j=2}^{n_1} \sum_{s=1}^{j-1} E(X_{1s}'X_{1j})^4 + \sum_{j=n_1+1}^{n_1+n_2} \sum_{s=1}^{n_1} E(X_{1s}'X_{2j-n_1})^4 \right.$$

$$\left. + \sum_{j=n_1+1}^{n_1+n_2} \sum_{s=n_1+1}^{j-1} E(X_{2s-n_1}'X_{2j-n_1})^4 \right\}$$

$$= O(n^{-8})(P_1 + P_2 + P_3),$$

where $P_1 = \sum_{j=2}^{n_1} \sum_{s=1}^{j-1} E(X_{1s}'X_{1j})^4$, $P_2 = \sum_{j=n_1+1}^{n_1+n_2} \sum_{s=1}^{n_1} E(X_{1s}'X_{2j-n_1})^4$ and

$$P_3 = \sum_{j=n_1+1}^{n_1+n_2} \sum_{s=n_1+1}^{j-1} E(X_{2s-n_1}'X_{2j-n_1})^4.$$

Let us consider $E(X_{1s}'X_{2j-n_1})^4$. Define $\Gamma_1'\Gamma_2 =: (v_{ij})_{m \times m}$ and note the following facts which will be used repeatedly in the rest of the Appendix:

$$\sum_{i,j=1}^{m} v_{ij}^4 \le \left( \sum_{i,j=1}^{m} v_{ij}^2 \right)^2 = \text{tr}^2(\Gamma_1'\Gamma_2\Gamma_2'\Gamma_1)$$

$$= \text{tr}^2(\Sigma_2\Sigma_1),$$

$$\sum_{i=1}^{m} \sum_{j_1 \ne j_2}^{m} (v_{ij_1}^2 v_{ij_2}^2) \le \left( \sum_{i,j=1}^{m} v_{ij}^2 \right)^2 = \text{tr}^2(\Sigma_2\Sigma_1),$$

$$\sum_{i_1 \ne i_2}^{m} \sum_{j_1 \ne j_2}^{m} v_{i_1j_1} v_{i_1j_2} v_{i_2j_1} v_{i_2j_2} \le \sum_{i_1 \ne i_2}^{m} v_{i_1i_2}^{(2)} v_{i_1i_2}^{(2)} \le \sum_{i_1,i_2=1}^{m} v_{i_1i_2}^{(2)} v_{i_1i_2}^{(2)},$$

$$\sum_{i_1,i_2=1}^{m} v_{i_1i_2}^{(2)} v_{i_1i_2}^{(2)} = \text{tr}(\Gamma_1'\Sigma_2\Gamma_1\Gamma_1'\Sigma_2\Gamma_1) = \sum_{i=1}^{m} v_{ii}^{(4)}$$

$$= \text{tr}\{(\Sigma_1\Sigma_2)^2\},$$

where $\Gamma_1'\Sigma_2\Gamma_1 = (v_{ij}^{(2)})_{m \times m}$ and $(\Gamma_1'\Sigma_2\Gamma_1)^2 = (v_{ij}^{(4)})_{m \times m}$.

From (3.1),

$$E(X_{1s}'X_{2j-n_1})^4 = \sum_{i=1}^{m} \sum_{j'=1}^{m} (3+\Delta)^2 v_{ij'}^4 + \sum_{i=1}^{m} (3+\Delta) \sum_{j_1 \ne j_2}^{m} v_{ij_1}^2 v_{ij_2}^2$$

$$+ \sum_{j'=1}^{m} (3+\Delta) \sum_{i_1 \ne i_2}^{m} v_{i_1j}^2 v_{i_2j}^2$$

$$+ 9 \sum_{i_1 \neq i_2}^{m} \sum_{j_1 \neq j_2}^{m} v_{i_1 j_1} v_{i_1 j_2} v_{i_2 j_1} v_{i_2 j_2}$$

$$= O\{\text{tr}^2(\Sigma_2 \Sigma_1)\} + O\{\text{tr}(\Sigma_2 \Sigma_1)^2\}.$$

Then we conclude

$$O(n^{-8})P_2 = \sum_{j=n_1+1}^{n_1+n_2} \sum_{s=1}^{n_1} [O\{\text{tr}^2(\Sigma_2 \Sigma_1)\} + O\{\text{tr}(\Sigma_2 \Sigma_1)^2\}]$$

$$= O(n^{-5})[O\{\text{tr}^2(\Sigma_2 \Sigma_1)\} + O\{\text{tr}(\Sigma_2 \Sigma_1)^2\}]$$

$$= o(\sigma_{n_1}^4).$$

We can also prove that $O(n^{-8})P_1 = o(\sigma_{n_1}^4)$ and $O(n^{-8})P_3 = o(\sigma_{n_1}^4)$ by going through a similar procedure. This completes the proof of the lemma. $\square$

PROOF OF THEOREM 1. We note equations (6.2) and (6.3) under conditions (3.4) and (3.5), respectively. Based on Corollary 3.1 of Hall and Heyde (1980), Lemmas 1, 2 and 3, it can be concluded that $T_{n_1}/\sigma_{n_1} \xrightarrow{d} N(0,1)$. This implies the desired asymptotic normality of $T_n$ under (3.4). Under (3.5), as $T_{n2}$ is the sum of two independent averages, its asymptotic normality can be attained by following the standard means. Hence the theorem is proved. $\square$

PROOF OF THEOREM 2. We only present the proof for the ratio consistency of $\widehat{\text{tr}(\Sigma_1^2)}$ as the proofs of the other two follow the same route. We want to show

(6.10) $\quad E\{\widehat{\text{tr}(\Sigma_1^2)}\} = \text{tr}(\Sigma_1^2)\{1 + o(1)\}$ and $\text{Var}\{\widehat{\text{tr}(\Sigma_1^2)}\} = o\{\text{tr}^2(\Sigma_1^2)\}.$

For notation simplicity, we denote $X_{1j}$ as $X_j$ and $\Sigma_1$ as $\Sigma$, since we are effectively in a one-sample situation.

Note that

$$\widehat{\text{tr}(\Sigma^2)} = \{n(n-1)\}^{-1}$$

$$\times \text{tr}\left[\sum_{j \neq k}^{n}\{(X_j - \mu)(X_j - \mu)'(X_k - \mu)(X_k - \mu)' \right.$$

$$- 2(\bar{X}_{(j,k)} - \mu)(X_j - \mu)'(X_k - \mu)(X_k - \mu)'\}$$

$$+ \sum_{j \neq k}^{n}\{2(X_j - \mu)\mu'(X_k - \mu)(X_k - \mu)'$$

$$- 2(\bar{X}_{(j,k)} - \mu)\mu'(X_k - \mu)(X_k - \mu)'\}$$

$$+ \sum_{j \neq k}^{n} \{(\bar{X}_{(j,k)} - \mu)(X_j - \mu)'(\bar{X}_{(j,k)} - \mu)(X_k - \mu)'\}$$

$$- \sum_{j \neq k}^{n} \{2(X_j - \mu)\mu'(\bar{X}_{(j,k)} - \mu)(X_k - \mu)'$$

$$- 2(\bar{X}_{(j,k)} - \mu)\mu'(\bar{X}_{(j,k)} - \mu)(X_k - \mu)'\}$$

$$+ \sum_{j \neq k}^{n} \{(X_j - \mu)\mu'(X_k - \mu)\mu' - 2(\bar{X}_{(j,k)} - \mu)\mu'(X_k - \mu)\mu'\}$$

$$+ \sum_{j \neq k}^{n} \{(\bar{X}_{(j,k)} - \mu)\mu'(\bar{X}_{(j,k)} - \mu)\mu'\}\Bigg]$$

$$=: \sum_{l=1}^{10} \text{tr}(A_l), \quad \text{say}.$$

It is easy to show that $E\{\text{tr}(A_1)\} = \text{tr}(\Sigma^2)$, $E\{\text{tr}(A_i)\} = 0$ for $i = 2, \ldots, 9$ and $E\{\text{tr}(A_{10})\} = \mu'\Sigma\mu/(n-2) = o\{\text{tr}(\Sigma^2)\}$. The last equation is based on (3.4). This leads to the first part of (6.10). Since $\text{tr}(A_{10})$ is nonnegative and $E\{\text{tr}(A_{10})\} = o\{\text{tr}(\Sigma^2)\}$, we have $\text{tr}(A_{10}) = o_p\{\text{tr}(\Sigma^2)\}$. However, to establish the orders of other terms, we need to derive $\text{Var}\{\text{tr}(A_i)\}$. We shall only show $\text{Var}\{\text{tr}(A_1)\}$ here. Derivations for other $\text{Var}\{\text{tr}(A_i)\}$ are similar.

Note that

$$\text{Var}\{\text{tr}(A_1)\} + \text{tr}^2(\Sigma^2)$$

$$= E\left[\frac{1}{n(n-1)} \text{tr}\left\{\sum_{j \neq k}^{n} (X_j - \mu)(X_j - \mu)'(X_k - \mu)(X_k - \mu)'\right\}\right]^2$$

$$= \frac{1}{n^2(n-1)^2} E\left[\text{tr}\left\{\sum_{j_1 \neq k_1}^{n} (X_{j_1} - \mu)(X_{j_1} - \mu)'(X_{k_1} - \mu)(X_{k_1} - \mu)'\right\}\right.$$

$$\left.\times \text{tr}\left\{\sum_{j_2 \neq k_2}^{n} (X_{j_2} - \mu)(X_{j_2} - \mu)'(X_{k_2} - \mu)(X_{k_2} - \mu)'\right\}\right].$$

It can be shown, by considering the possible combinations of the subscripts $j_1, k_1, j_2$ and $k_2$, that

$$\text{Var}\{\text{tr}(A_1)\} = 2\{n(n-1)\}^{-1} E\{(X_1 - \mu)'(X_1 - \mu)\}^4$$

(6.11)
$$+ \frac{4(n-2)}{n(n-1)} E\{(X_1 - \mu)'\Sigma(X_1 - \mu)\}^2 + o\{\text{tr}^2(\Sigma^2)\}$$

$$=: \frac{2}{n(n-1)}B_{11} + \frac{4(n-2)}{n(n-1)}B_{12} + o\{\operatorname{tr}^2(\Sigma^2)\},$$

where

$$B_{11} = E(Z_1'\Gamma'\Gamma Z_2)^4 = E\left(\sum_{s,t=1}^m z_{1s}\nu_{st}z_{2t}\right)^4$$

$$= E\left(\sum_{s_1,s_2,s_3,s_4,t_1,t_2,t_3,t_4=1}^m \nu_{s_1t_1}\nu_{s_2t_2}\nu_{s_3t_3}\nu_{s_4t_4} z_{1s_1}z_{1s_2}z_{1s_3}z_{1s_4}z_{2t_1}z_{2t_2}z_{2t_3}z_{2t_4}\right)$$

and

$$B_{12} = E(Z_1'\Gamma'\Gamma\Gamma'\Gamma Z_1)^2 = E\left(\sum_{s,t=1}^m z_{1s}u_{st}z_{1t}\right)^2$$

$$= E\left(\sum_{s_1,s_2,t_1,t_2=1}^m u_{s_1t_1}u_{s_2t_2}z_{1s_1}z_{1s_2}z_{1t_1}z_{1t_2}\right).$$

Here $\nu_{st}$ and $u_{st}$ are, respectively, the $(s,t)$ element of $\Gamma'\Gamma$ and $\Gamma'\Sigma\Gamma$.

Since $\operatorname{tr}^2(\Sigma^2) = (\sum_{s,t=1}^m \nu_{st}^2)^2 = \sum_{s_1,s_2,t_1,t_2=1}^m \nu_{s_1t_1}^2 \nu_{s_2t_2}^2$ and $\operatorname{tr}(\Sigma^4) = \sum_{t_1,t_2=1}^m u_{t_1t_2}^2$. It can be shown that $A_{11} \leq c\operatorname{tr}^2(\Sigma^2)$ for a finite positive number $c$ and hence $\{n(n-1)\}^{-1}B_{11} = o\{\operatorname{tr}^2(\Sigma^2)\}$. It may also be shown that

$$B_{12} = 2\sum_{s,t=1}^m u_{st}^2 + \sum_{s,t=1}^m u_{ss}u_{tt} + \Delta\sum_{s=1}^m u_{ss}^2$$

$$= 2\operatorname{tr}(\Sigma^4) + \operatorname{tr}^2(\Sigma^2) + \Delta\sum_{s=1}^m u_{ss}^2$$

$$\leq (2+\Delta)\operatorname{tr}(\Sigma^4) + \operatorname{tr}^2(\Sigma^2).$$

Therefore, from (6.11),

$$\operatorname{Var}\{\operatorname{tr}(A_1)\} \leq \frac{2}{n(n-1)}c\operatorname{tr}^2(\Sigma^2) + \frac{4(n-2)}{n(n-1)}\{(2+\Delta)\operatorname{tr}(\Sigma^4) + \operatorname{tr}^2(\Sigma^2)\}$$

$$= o\{\operatorname{tr}^2(\Sigma^2)\}.$$

This completes the proof. □

**Acknowledgments.** We are grateful to two reviewers for valuable comments and suggestions which have improved the presentation of the paper. We also thank Dan Nettleton and Peng Liu for useful discussions.

## REFERENCES

ANDERSON, T. W. (2003). *An Introduction to Multivariate Statistical Analysis*. Wiley, Hoboken, NJ. MR1990662

ABRAMOVICH, F., BENJAMINI, Y., DONOHO, D. L. and JOHNSTONE, I. M. (2006). Adaptive to unknown sparsity in controlling the false discovery rate. *Ann. Statist.* **34** 584–653. MR2281879

BAI, Z. and SARANADASA, H. (1996). Effect of high dimension: By an example of a two sample problem. *Statist. Sinica* **6** 311–329. MR1399305

BARRY, W., NOBEL, A. and WRIGHT, F. (2005). Significance analysis of functional categories in gene expression studies: A structured permutation approach. *Bioinformatics* **21** 1943–1949.

BENJAMINI, Y. and HOCHBERG, Y. (1995). Controlling the false discovery rate: A practical and powerful approach to multiple testing. *J. Roy. Statist. Soc. Ser. B* **57** 289–300. MR1325392

BENJAMINI, Y. and YEKUTIELI, D. (2001). The control of the false discovery rate in multiple testing under dependency. *Ann. Statist.* **29** 1165–1188. MR1869245

CHEN, S. X. and QIN, Y.-L. (2008). *A Two Sample Test for High Dimensional Data with Applications to Gene-Set Testing*. Research report, Dept. Statistics, Iowa State Univ.

CHIARETTI, S., LI, X. C., GENTLEMAN, R., VITALE, A., VIGNETTI, M., MANDELLI, F., RITZ, J. and FOA, R. (2004). Gene expression profile of adult T-cell acute lymphocytic leukemia identifies distinct subsets of patients with different response to therapy and survival. *Blood* **103** 2771–2778.

DUDOIT, S., KELES, S. and VAN DER LAAN, M. (2008). Multiple tests of association with biological annotation metadata. *Inst. Math. Statist. Collections* **2** 153–218.

EFRON, B. and TIBSHIRANI, R. (2007). On testing the significance of sets of genes. *Ann. Appl. Stat.* **1** 107–129. MR2393843

FAN, J., HALL, P. and YAO, Q. (2007). To how many simultaneous hypothesis tests can normal, Student's t or bootstrap calibration be applied. *J. Amer. Statist. Assoc.* **102** 1282–1288. MR2372536

FAN, J., PENG, H. and HUANG, T. (2005). Semilinear high-dimensional model for normalization of microarray data: A theoretical analysis and partial consistency. *J. Amer. Statist. Assoc.* **100** 781–796. MR2201010

GENTLEMAN, R., IRIZARRY, R. A., CAREY, V. J., DUDOIT, S. and HUBER, W. (2005). *Bioinformatics and Computational Biology Solutions Using R and Bioconductor*. Springer, New York. MR2201836

HALL, P. and HEYDE, C. (1980). *Martingale Limit Theory and Applications*. Academic Press, New York. MR0624435

HUANG, J., WANG, D. and ZHANG, C. (2005). A two-way semilinear model for normalization and analysis of cDNA microarray data. *J. Amer. Statist. Assoc.* **100** 814–829. MR2201011

KOSOROK, M. and MA, S. (2007). Marginal asymptotics for the "large $p$, small $n$" paradigm: With applications to microarray data. *Ann. Statist.* **35** 1456–1486. MR2351093

LEDOIT, O. and WOLF, M. (2002). Some hypothesis tests for the covariance matrix when the dimension is large compare to the sample size. *Ann. Statist.* **30** 1081–1102. MR1926169

NEWTON, M., QUINTANA, F., DEN BOON, J., SENGUPTA, S. and AHLQUIST, P. (2007). Random-set methods identify distinct aspects of the enrichment signal in gene-set analysis. *Ann. Appl. Stat.* **1** 85–106. MR2393842

Portnoy, S. (1986). On the central limit theorem in $R^p$ when $p \to \infty$. *Probab. Theory Related Fields* **73** 571–583. MR0863546

Recknor, J., Nettleton, D. and Reecy, J. (2008). Identification of differentially expressed gene categories in microarray studies using nonparametric multivariate analysis. *Bioinformatics* **24** 192–201.

Schott, J. R. (2005). Testing for complete independence in high dimensions. *Biometrika* **92** 951–956. MR2234197

Storey, J., Taylor, J. and Siegmund, D. (2004). Strong control, conservative point estimation and simultaneous conservative consistency of false discovery rates: A unified approach. *J. R. Stat. Soc. Ser. B Stat. Methodol.* **66** 187–205. MR2035766

Tracy, C. and Widom, H. (1996). On orthogonal and symplectic matrix ensembles. *Comm. Math. Phys.* **177** 727–754. MR1385083

Van der Laan, M. and Bryan, J. (2001). Gene expression analysis with the parametric bootstrap. *Biostatistics* **2** 445–461.

Yin, Y., Bai, Z. and Krishnaiah, P. R. (1988). On the limit of the largest eigenvalue of the large-dimensional sample covariance matrix. *Probab. Theory Related Fields* **78** 509–521. MR0950344

Department of Statistics
Iowa State University
Ames, Iowa 50011-1210
USA
and
Guanghua School of Management
Peking University
Beijing 100871
China
E-mail: songchen@iastate.edu

Department of Statistics
Iowa State University
Ames, Iowa 50011-1210
USA
E-mail: qinyl@iastate.edu